\documentclass[a4paper,12pt]{article}

\usepackage{float}
\usepackage{times}
\usepackage{graphicx}

\renewcommand\refname{References and Notes}
\renewcommand{\figurename}{\textbf{Figure}}
\usepackage[labelfont=bf]{caption}
\DeclareCaptionLabelSeparator{pipe}{ $\mathbf{\vert}$ }
\title{Exploring the quantum speed limit with computer games}
\author{Jens Jakob W. H. S\o rensen, Mads Kock Pedersen, Michael Munch, Pinja Haikka,\\
Jesper Halkj{\ae}r Jensen, Tilo Planke, Morten Ginnerup Andreasen,\\
Miroslav Gajdacz, Klaus M\o lmer, Andreas Lieberoth,\\
Jacob F. Sherson$^\dagger$, and Quantum Moves players\\
\\
\normalsize{Department of Physics and Astronomy, Aarhus University, Denmark.}\\
\\
\normalsize{$^\dagger$sherson@phys.au.dk}}
\date{\today}

\begin{document}
\baselineskip12pt
\maketitle

\begin{abstract}
Humans routinely solve problems of immense computational complexity by intuitively
forming simple, low-dimensional heuristic strategies \cite{McLeod1996,Gigerenzer1999}. Citizen science exploits this
ability by presenting scientific research problems to non-experts. Gamification is an effective
tool for attracting citizen scientists to provide solutions to research problems. While citizen science games Foldit \cite{Cooper2010}, EteRNA \cite{Lee2014} and EyeWire \cite{Kim2014} have been used successfully to study protein and RNA folding and neuron mapping, so far gamification has not been applied to problems in quantum physics. Does the fact that everyday experiences are based on classical physics hinder the use of non-expert citizen scientists in the realm of quantum mechanics? Here we report on Quantum Moves, an online platform gamifying optimization problems in quantum physics. We show that  human players are able to find solutions to difficult problems associated with the task of quantum computing \cite{Weitenberg2011a}. Players succeed where purely numerical optimization fails, and analyses of their solutions provide insights into the problem of optimization of a more profound and general nature. Based on player strategies, we have thus developed a new, few-parameter heuristic optimization method which efficiently outperforms the most prominent established numerical methods. The numerical complexity associated with time-optimal solutions increases for shorter process durations. To better understand this, we have made a low-dimensional rendering of the optimization landscape. These studies show why traditional optimization methods fail near the quantum speed limit \cite{Mandelstam1945,Caneva2011,Gajdacz2014}, and they bring promise that combined analyses of optimization landscapes and heuristic solution strategies may benefit wider classes of optimization problems in quantum physics and beyond. 
\end{abstract}

Quantum physics holds the potential of unprecedented technological advances in the realms of computing \cite{Monroe2002} and simulations \cite{Lewenstein2007}. To ensure functionality, all quantum operations must be executed near perfection, requiring highly optimized operations with fidelities above $F\geq 0.999$ \cite{Devitt2013}. Given the high dimensionality of quantum optimization problems, one could expect these \textit{quantum optimal control} problems to be impractically difficult to solve. However, assuming full controllability, quantum optimization problems are benign since all local maxima are also global maxima \cite{Rabitz}. Tailored local optimizers, such as the gradient-based Krotov Algorithm (KA), solve these problems \cite{Sklarz2002}.

Quantum computing operations must be executed faster than typical decoherence times to ensure functionality. However, there is a shortest process duration with perfect fidelity, denoted the quantum speed limit (QSL) \cite{Mandelstam1945}, which imposes a fundamental limit on the process duration and hence quantum computation. With limited duration, the quantum optimization problem loses the favourable properties stated earlier and local optimizers are not almost guaranteed to converge. In this case the prevalent and hitherto successful paradigm is using multistart of such local algorithms \cite{Ugray2007,Caneva2009,Caneva2011}. For a time-dependent Hamiltonians the QSL is conventionally computed by assuming that it coincides with the process duration for which the multistarting of local optimization fails \cite{Caneva2009, Caneva2011}. We demonstrate that this assumption is not necessarily true.

High-dimensional optimization problems are often solved by humans using simple heuristic strategies. Common examples are visual pattern recognition \cite{Bilalic2010} and catching a flying ball subject to wind and air resistance \cite{McLeod1996}. Citizen science projects such as Foldit \cite{Cooper2010}, EyeWire \cite{Kim2014} and Galaxy Zoo \cite{Lintott2011} employ these human skills to solve highly complex research problems via gamification. We ask whether citizen science projects can be extended from these puzzle and pattern recognition tasks to dynamic challenges, and whether this approach can be implemented on quantum physics problems. In order to investigate these questions we created Quantum Moves\footnote{Found at www.scienceathome.org}, which presents quantum computing operations as games. 

The quantum computer architecture we gamify employs neutral atoms trapped in optical lattices \cite{Lewenstein2007}. By driving the superfluid-Mott insulator transition, samples of hundreds of atoms have been trapped in optical lattices with an unprecedented purity \cite{Bakr2010,Sherson2010}. Their regular spacing offers the possibility of creating a scalable quantum computer, which is a challenge for other available architectures \cite{Ladd2010}. Several proposals for the implementation of quantum computing in this system have been proposed, mainly using long range gates or contact interactions \cite{Lewenstein2007}. Here we investigate an architecture \cite{Weitenberg2011a} where contact interaction is achieved by moving atoms on top of each other using a so-called optical tweezer \cite{Weitenberg2011b,Kaufman2015}, a tightly focused off-resonant laser beam. Finding the optimal motion of the tweezer from one lattice site to another is a difficult problem when the available time is close to the QSL, and it is this transfer problem that is introduced to the players of Quantum Moves through different challenges.

\begin{figure}[H]
\centering
\includegraphics[width=0.5\textwidth]{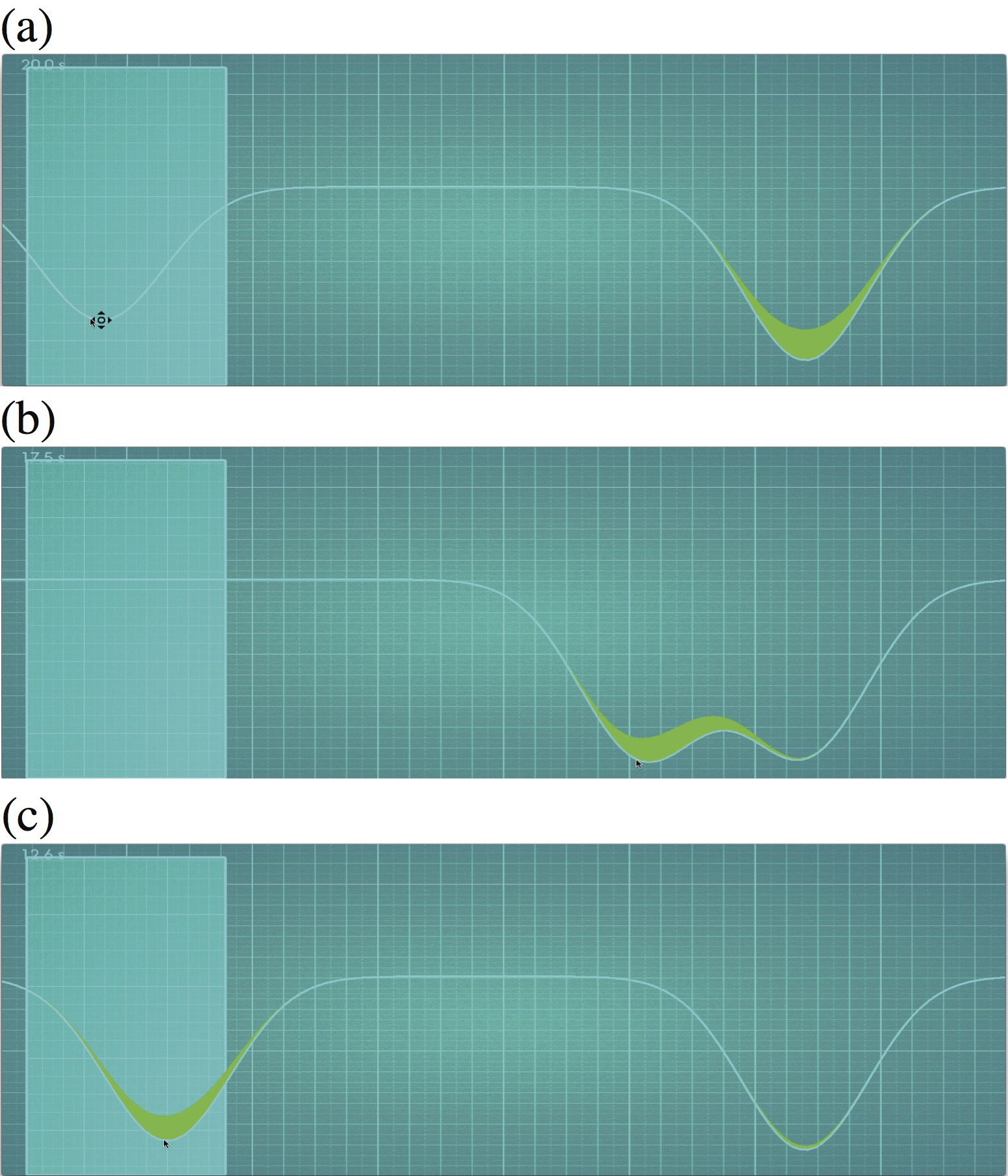}
\caption{\textbf{The BringHomeWater (BHW) challenge as seen by the player.} The atom is represented by the square of its wave function, $|\psi(x,t)|^2$, shown as a green liquid. The blue curve represents the potential felt by the atom. Panel \textbf{b} The controllable tweezer is initially on the left and the atom is trapped in the right static potential. The player controls the optical tweezer by moving a computer cursor, picks up the atom (\textbf{b}), and drags it back to the target area (\textbf{c}), marked by a cyan rectangle, to collect points.}
\label{fig:fig1}
\end{figure}

Here we present the results of one challenge called BringHomeWater (BHW). The aim of BHW is to move the optical tweezer into a region where an atom is trapped in a fixed potential well, collect the atom and move it back to a target area as quickly as possible. Fast movement introduces excitations in the state of the atom, which is described by a quantum mechanical wave function $\psi(x,t)$ and visualized as a sloshy liquid---see Fig. \ref{fig:fig1} for the player view of BHW and Methods for a description of the game interface and player demographics. The excitations must be stabilized before reaching the target area to attain high fidelity of the underlying quantum process.

BHW belongs to the class of one-dimensional quantum optimal control of individual atoms with one or two control parameters. This class of problems is of great interest, has been studied extensively within the past decade (see e.g. Refs. \cite{Weitenberg2011a,Gajdacz2014,deChiara2008}) and is highly relevant for the experimental platforms in Refs. \cite{Sherson2010,jager2014}. We believe our gamification strategy is extendable to these problems and other physical interactions including, but not limited to, many-body dynamics in optical lattices \cite{doria2011} and Bose-Einstein condensates \cite{jager2014}.

Recall that the standard approach for solving a problem such as BHW is to use a multistart of tailored local optimization algorithms, such as the Krotov algorithm (KA) \cite{Sklarz2002,Caneva2009,Caneva2011}. The choice of an initial seed is a central challenge in such complex optimizations. In BHW the high dimensional optimization space can only be searched sparsely, and we used multiple methods for creating initial seeds for the KA. The most successful method KASS employs linear combinations of sinusoidal functions as seed trajectories (see Methods). Using the KA on these seeds, high fidelity solutions are readily found for process durations longer than $T=0.40$ (units defined in Methods). Iteratively, solutions with shorter durations are found by contracting the solutions with a slightly longer duration and using them as seeds for the KA. This procedure is called a sweep, and it traces out entire families of solutions. The best results of KASS are displayed in red in Fig. \ref{fig:fig2}a.  This method locates an ostensible QSL at  $T_{\mathrm{QSL}}^{\mathrm{num}}=0.29$ after approximately $7.4 \cdot 10^8$ trials. This KASS optimization formed our most successful bare computer optimization method, outperforming also the acclaimed CRAB algorithm \cite{Caneva2011b}.

KASS and CRAB fit the prevalent paradigm of multistarting of local optimizers, while global optimizers are rarely used in quantum optimal control. To investigate the BHW problem using a global optimizer we chose the so-called differential evolution algorithm (DE) due to its demonstrated success in quantum problems \cite{Zahedinejad2014}. For duration $T=0.40$, DE performed worse than the Krotov optimization with identical computational resources. The poor performance of DE can be attributed to the very high dimensionality of the optimization space, and to the scarceness of good solutions found therein.

Players were trained in a series of introductory levels before reaching the scientific challenges. Training levels equip the players with a range of skills; successful solutions of the scientific challenges require a holistic combination of these skills. In total, all Quantum Moves games have been played about 500,000 times by roughly 10,000 players. BHW is the most played scientific level (approx. 12,000 plays by 300 players) and is therefore analysed in detail here. Player results span many different fidelities and process durations---see the dots in Fig. \ref{fig:fig2}a. Remarkably, players trace out a region in the vicinity of $T_{\mathrm{QSL}}^{\mathrm{num}}$ very similar to that of the best obtained numerical results, despite the fact that the numerical optimization used roughly 100,000 times more trials than the players. For short durations, players find even \textit{better} solutions than the numerical optimization, albeit with imperfect fidelities.

This result inspired us to introduce a powerful hybrid Computer-Human Optimization (CHOP) scheme, in which we use the players' intuitive solutions as seeds for the local numerical optimization. CHOP was applied to the top 70 \% of the player solutions with durations shorter than $T=0.40$. The results of these player-seeded optimization sweeps are shown in Fig. \ref{fig:fig2}b. It turns out that the CHOP solutions bunch in groups, which we denote as clans and discuss in detail later. In Fig. \ref{fig:fig2}a we illustrate the best solution families of the two dominant clans (blue and yellow curves). The QSL found by CHOP is $T_{\mathrm{QSL}}^{\mathrm{num}}=0.20$, a vast improvement of the value obtained from the bare computer optimization using multistarted local optimization. This result clearly demonstrates that BHW is a quantum control problem without the benign properties discussed previously and questions the common assumption that the QSL can be found numerically as the duration for which local optimization fails \cite{Caneva2009}.

To better understand the ability of the players and CHOP to give better solutions than KASS, we examine the clustering of the CHOP solutions. We introduce a measure describing the distance between two solutions for a particular process duration \textit{T}:
\begin{equation}
D_{j,k}=\frac{1}{T}\int_0^T\langle f_{jk}|f_{jk}\rangle \mathrm{d}t,
\label{distance}
\end{equation}
where $|f_{jk}\rangle=|\psi_j(x,t)\rangle -|\psi_k(x,t)\rangle$ is the difference between two wave functions evaluated along the path. Figure \ref{fig:fig3}a shows a distance map of the CHOP solutions. Two dominant clans stand out as distinct regions in the distance map---See Methods for details on the identification of clans. The solutions of the BHW transport problem corresponding to these two clans are shown in Fig. \ref{fig:fig3}b. Using the first solution, marked in yellow in Fig. \ref{fig:fig3}b, the atom is collected by \textit{tunneling} the wave function into a tweezer potential placed on the left hand side of the static potential. In the second class of \textit{shoveling} solutions, marked in blue in Fig. \ref{fig:fig3}b, the tweezer is moved past the position of the atom so that the overlap of the tweezer and the static potential forms a strong potential gradient accelerating the atom towards the target. Fast non-adiabatic solutions must spread the wave function into different energy eigenstates. In Fig.~\ref{fig:fig3}c we see that the solutions in both clans populate many different eigenstates. The tweezer must then be shaken periodically when approaching the target area to bring the atom back to the desired ground state---see Methods for more details.  \textit{A priori}, it was not obvious that two strategies, corresponding to distinctly different physical phenomena, should exist for this kind of a problem. It is also worth stressing that players explore both solutions despite having no or little prior knowledge of quantum mechanical phenomena such as tunneling.

The difficulty of a particular optimization problem can often be assessed with knowledge of the topology or ruggedness of the so-called optimization or fitness landscape. This has been established for problems in quantum optimal control \cite{Roslund2009} and in other fields \cite{Oana2014}. In this landscape the quality of the solution for each set of the control parameters is represented as the height. If many global optima are distributed across the landscape, local gradient-based search is often sufficient to find the highest peak. On the other hand, if the landscape is very rugged and contains many local maxima, local methods will in general fail. In order to understand our problem in this terminology we did a dimensional reduction of the data from the player seeded solutions and the bare numerical optimization. We assign points to individual solutions in a two-dimensional space such that the Euclidean distance between two points approximates the Manhattan type distance between two solutions in the full high dimensional optimization space (see Methods). The height of the landscape quantifies the fidelity of a solution.

For process durations $T=0.40$ and $T\simeq 0.17$, above and below $T_{\mathrm{QSL}}^{\mathrm{num}}$, the landscapes are visualized in Fig. \ref{fig:fig4}a and \ref{fig:fig4}b, respectively. Green areas in the landscape mark the player-generated CHOP solutions. For the longer duration (Fig. \ref{fig:fig4}a), global maxima indeed spread across the optimization landscape, explaining the success of Krotov-based methods. For the shorter duration (Fig. \ref{fig:fig4}b) all high-fidelity solutions lie in the green CHOP region. This explains the failure of the Krotov-based methods and difficulty in locating the true $T_{\mathrm{QSL}}$, since the global optima are no longer spread across the landscape. Thus Fig. \ref{fig:fig4}a and \ref{fig:fig4}b demonstrate the dramatic change in the landscape as the duration is decreased. It is worth stressing that the CHOP region is tiny in the high-dimensional landscape and therefore easily missed even by elaborate seeding strategies. CHOP outperforms the numerical optimization because players are able to heuristically identify the regions of high fidelity in the high-dimensional optimization landscape, thereby finding the best seeds for Krotov-type local algorithms.

One long term vision of our work is to circumvent the need for gamification by learning how players form their successful low-dimensional heuristic strategies, and to incorporate this into autonomous optimization algorithms. As a first step in this direction, we introduce here a Heuristically Initialized Local Optimizer (HILO) algorithm. HILO parametrizes the player solutions in a low-dimensional subspace, whilst retaining the main features of good seeds. A local search algorithm can then move beyond this subspace to find optimal solutions. More specifically, inspired by the player solutions we constructed a three-dimensional parameterization consisting of moving the tweezer i) right, ii) slowly left, and iii) quickly left (see Methods for details). Paths from this three-dimensional space for $T=0.15$ were used as seeds for the KA and iteratively applied to shorter and longer durations with a sweep (purple curve in Fig. \ref{fig:fig2}a). Seeds were taken from the low-dimensional subspace using a simple direct search. HILO finds the lowest QSL at $T_{\mathrm{QSL}}^{\mathrm{num}} = 0.19$, outperforming even the best CHOP strategies. The solutions from HILO are shown as a landscape in Fig.~\ref{fig:fig4}c. Initialized in a low-dimensional seed space, HILO efficiently explores a smaller but more optimal volume of the global optimization space. Only parametrizations which accurately capture the nature of efficient solution strategies at short durations will lead to efficient optimization. 
The intuition gained from the players was pivotal as the parametrization used in HILO emerged from the CHOP solutions and the physical interpretation of the player-based solution strategies. This makes CHOP a crucial precursor to HILO.
We stress that any optimal control effort requires a substantial amount of optimizations of different seeds to find good solutions. After this initial work, HILO could be readily implemented by first applying our clustering methods to identify potential clans. If successful, analysis of the clans would allow for a formulation of the low-dimensional parametrization used in HILO.

Finally, we discuss another central component of quantum optimal control theory, namely the dependence of the fidelity on the process duration around the QSL. This has previously been associated with a universal-type $\sin^2\bigl(\frac{\pi}{2} T/T_{\mathrm{QSL}}\bigr)$ behavior \cite{Caneva2011}. Here we do not observe this behavior (see Extended Data Fig. \ref{fig:figS1}a and Methods for details). We attribute this to variations in the so-called direct Hilbert speed \cite{Gajdacz2014}. The variations in our problem may arise from the sequential nature of the good strategies, such as the three steps in the HILO-parametrization, where each sequence has a different scaling of fidelity with time.

Using the gamified interface in Quantum Moves, players with little or no training in quantum physics not only provided high quality solutions, but also enabled the extraction of the underlying physical strategies. This success encourages the pursuit of other quantum research games, as well as \textit{dynamic} games in other fields. Our future work will focus on both extending the gamification strategy and also the more general classification methodology to new types of control problems such as those presented in Refs. \cite{doria2011}, \cite{jager2014} and \cite{deChiara2008}. For all these different physical interactions we expect duration-constrained problems to exhibit complex optimization landscapes. Here players are expected to provide a global overview beyond typical local optimization. Another interesting topic of future research will be the exploration of how the players form strategies. Analyses of the player data will enable us to identify important features of quantum optimal control problems and allow an efficient dimensionality reduction. This will form the first steps in the efficient training of modern machine learning algorithms. As a first step in this direction, we have collaborated with cognitive researches to create a new game called Quantum Minds also available at www.scienceathome.org.

\bibliographystyle{naturemag}

\section*{End Notes}
\textbf{Acknowledgments}
We thank J. Rafner for graphical support and J. Jarecki, O. Vuculescu and C. Bergenholtz for fruitful discussions.
This work was supported by the European Research Council, Lundbeck Foundation, Aarhus University Research Foundation, Templeton Foundation, The Danish Council for Independent Research, Villum Foundation and Carlsberg Foundation.

\textbf{Author Contributions} 
All authors contributed to the construction of the online game platform and the effort to enlist users. J.J.W.H.S., M.K.P., T.P., M.G.A., M.G., K.M, and J.F.S. participated in the numerical analysis of the player and computer results. All authors contributed to the writing of the manuscript.

\textbf{Author Information} Reprints and permissions information is available at
www.nature.com/reprints. The authors declare no competing financial interests. Correspondence and requests for materials should be
addressed to J.F.S. (sherson@phys.au.dk).

\section*{Figures}

\begin{figure}[H]
\centering
\includegraphics[width=0.5\textwidth]{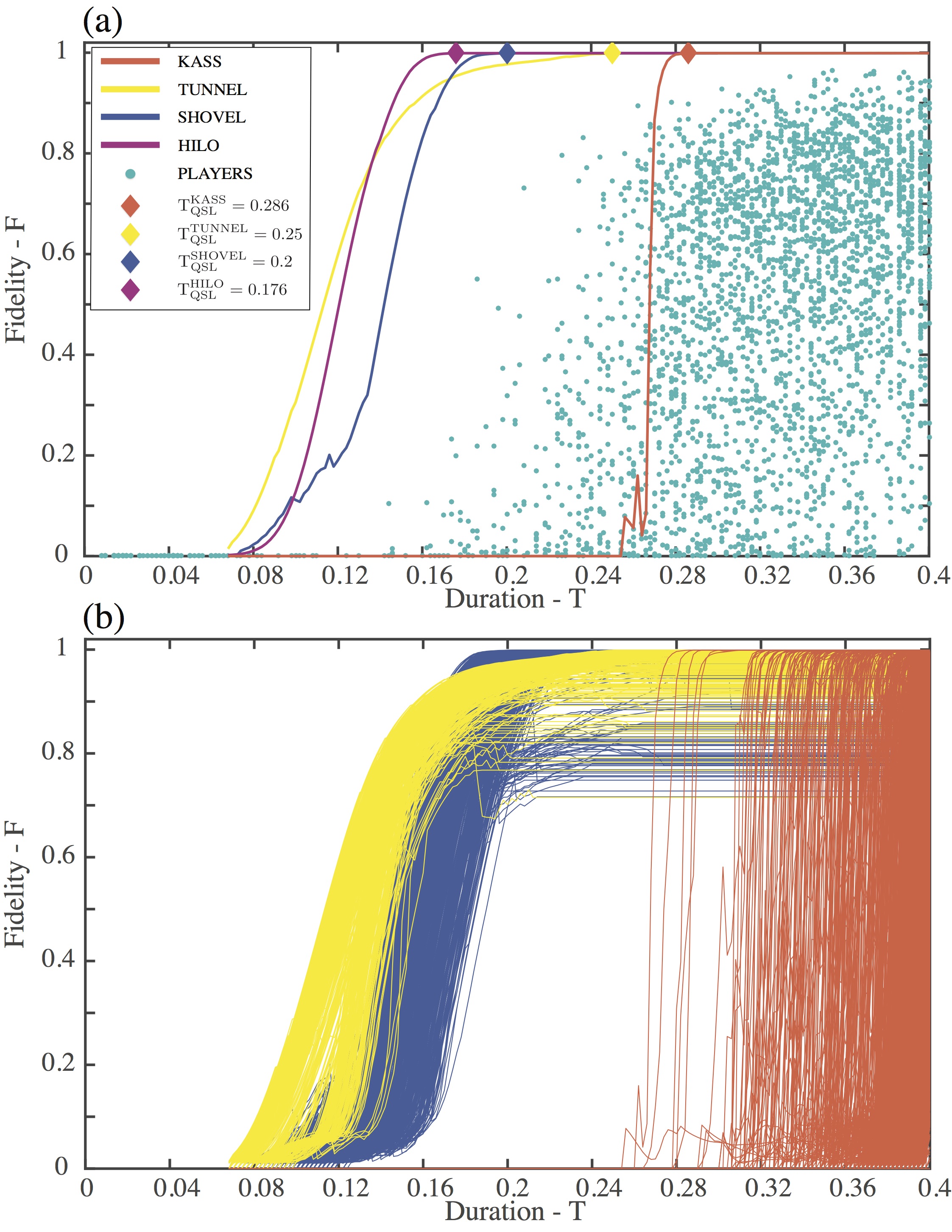}
\caption{\textbf{Fidelities of the transport problem for different solution durations.} Panel \textbf{a} shows a subset of the solutions found by the players as turquoise dots. The curves show the best solutions found by computer optimizations for each duration. The optimizations shown are KASS (red), shovelling (blue), tunneling (yellow), and HILO (purple). The diamonds mark the shortest duration with optimal fidelity ($F \geq 0.999$) for each optimization method. Panel \textbf{b} shows the sweeps from seeds that are generated by players (yellow and blue) and computers (red). Player solutions divide into shovelling (blue) and tunneling (yellow) clans.}
\label{fig:fig2}
\end{figure}

\begin{figure}[H]
\centering
\includegraphics[width=\textwidth]{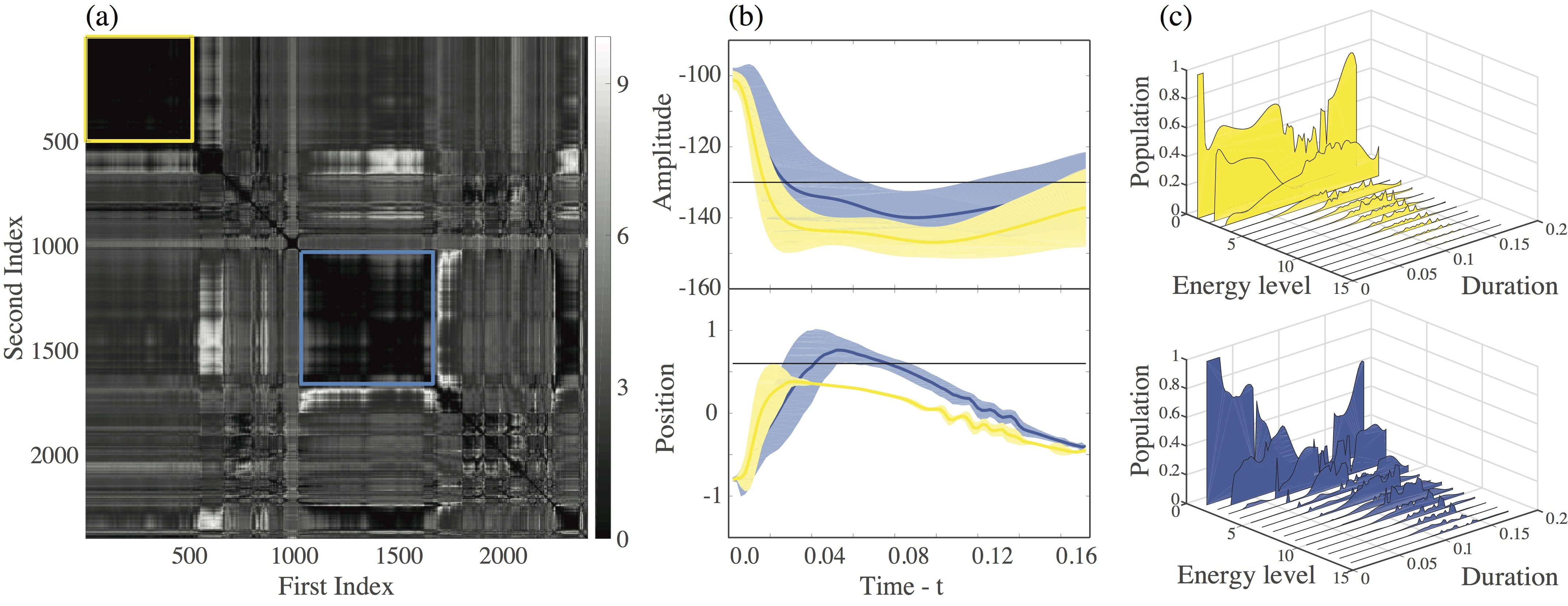}
\caption{\textbf{Shovelling and tunneling clans.} The shovelling (tunneling) clan is marked blue (yellow) throughout. Panel (a) is a distance map, showing the distances between CHOP solutions as defined by Eq. (\ref{distance}). Boundaries of the clans are marked with coloured squares. Panel \textbf{b} shows the average trajectories followed by the clans as thick lines, and a single standard deviations thereof as translucent areas. Trajectories are divided into the tweezer amplitude (top) and tweezer position (bottom). Panel \textbf{b} shows the wave function fraction (population) in the different instantaneous energy eigenstates.}
\label{fig:fig3}
\end{figure}

\begin{figure}[H]
\centering
\includegraphics[width=\textwidth]{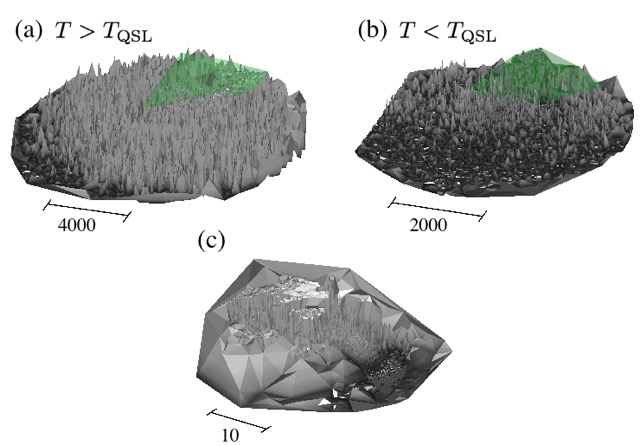}
\caption{\textbf{Optimization landscapes}. Panels \textbf{a} and \textbf{b} show the 2D rendering of the high-dimensional optimization landscape for process durations $T = 0.40$ and $T = 0.17$, respectively. Green areas mark the space probed by CHOP solutions. Panel \textbf{b} is the low-dimensional HILO landscape.}
\label{fig:fig4}
\end{figure}

\newpage
\section*{Methods}
\subsection*{Quantum Moves game interface}
In the online computer game Quantum Moves, the potential created by the optical tweezer is represented along the horizontal (position) and vertical (energy) axes (see Fig. \ref{fig:fig2}a). The potential is updated as the player uses the computer cursor to control the position and intensity of the optical tweezer. The sloshy liquid seen by the player on top of the tweezer potential is the probability density of the atom position. The probability density ($|\psi(x,t)|^2$) is defined by the quantum mechanical wave function $\psi(x,t)$. The wave function is updated using the time-dependent Schr\" odinger equation, according to the player-controlled potential. Thus, the time evolution is both observed and affected in real time by the player. The time is scaled approximately by a factor of $3 \cdot 10^{4}$, so that atomic evolution at microsecond time scale is experienced over a typical 10 second execution time of the individual game. 

The aim of the games is to find the fastest possible path that transfers the atom into the (stationary) ground state with unit fidelity, $F=1$, in the specified target region. When the tweezer is moved quickly, the atomic probability distribution begins to visibly slosh in a fluid-like fashion and the players may readily adopt strategies to prevent or control excitation (sloshing) of the wave. To encourage players to search for fast solutions necessary for realistic quantum computations, and to probe the region around the QSL, we introduce a time penalty in the game structure. For a thorough discussion on the Quantum Moves player demographics, including their scientific backgrounds, we refer the reader to \cite{Lieberoth2014} and for a similar analysis of Zooniverse.org players see \cite{Sauermann2015}.

\subsection*{KASS \label{sec:crb}}
We tried different seeding and optimization strategies, and the most successful was a Krotov algorithm using sinusoidal seed functions and sweeps over the total duration (KASS). KASS applies a sweep to an initial random seed $u_i(t)$, where $u_1(t)=\mathcal{A}(t)$ is the tweezer amplitude and $u_2(t)=x(t)$ is the tweezer position. The initial random seed is
\begin{equation}
u_i(k \cdot \delta t)=w_i(k \cdot \delta t)+\sum_{n=0}^{N-1} X [n] \sin\biggl(\frac{n \pi k}{N}\biggr), \label{sinbasis}
\end{equation}
where \textit{N} is the number of time slices in the path and $\delta t=0.002$ is the time discretization. $w_i$ is the motion at constant speed to the initial position of the atom and back again for $\mathcal{A}(t)=-100$. The amplitudes $X[n]$ were selected with a random sign and a norm decreasing as a function of the frequency. We generated roughly $2400$ seeds using the random summation series and apply the Krotov algorithm to these. The next step in KASS is to apply a sweep to the random seeds $u_i(t)$:
\begin{enumerate}
\item The optimized solution for the tweezer depth and position is linearly contracted in time by $\delta t$, $\mathcal{A}(t), x(t)\rightarrow \mathcal{A}(a\cdot t),x(a\cdot t)$ for $0<a<1$.
\item The Krotov algorithm is applied using the contracted solution as seed.
\item Repeat 1-2 until a minimal duration of $T=0.07$ is reached.
\end{enumerate}
Each of the $2400$ optimizations took 6 hours of calculation individually.

\subsection*{Schr\"odinger equation in dimensionless units \label{sec:units}}
The time-dependent Schr\" odinger equation solved in Quantum Moves has the dimensionless form,
\[
-\frac{1}{2}\frac{\mathrm{d}^2 \psi}{\mathrm{d}x^2}+V(x,t)\psi=i \frac{\mathrm{d} \psi}{\mathrm{d} t},
\]
where the position $x$ is measured in units of  $l_u=\lambda/2=532 \mathrm{nm}$. This is a typical period of optical lattice potentials used in experiments with ${}^{87}\mathrm{Rb}$ atoms \cite{Sherson2010}. The potential energy $V(x,t)$ is given in units of $E_u=1.70  \mathrm{peV}$, about $1/5$ of the so-called recoil energy in the lattice potential. Time is given in units of $t_u=\hbar/E_u = 0.39 \mathrm{ms}$. In these dimensionless units, our optical tweezer is parametrized as
\[
V_{\mathrm{tweezer}}=\mathcal{A}\exp\Biggl(-\frac{2.0(x-x_0)^2}{w_0^2}\Biggr),
\]
where $x_0$ is the position and $w_0^2$ is the waist. In the BHW challenge we have chosen $\mathcal{A}=130$ for the depth of the static tweezer potential and $w_0=0.25$ for the width of both tweezer potentials. The time-dependent Schr\" odinger equation is solved using the split-step method.

\subsection*{Identification of solution clans}
The distance between solutions was calculated using integral of Eq. (\ref{distance}) for  $T\simeq 0.17$, the duration for which the two best CHOP families in Fig. \ref{fig:fig2}b intersect.  The solutions were then sorted by picking one state at random, and choosing the next state on the list to be the closest one in the sense of Eq. (\ref{distance}). This process is iterated with the $n+1$-st solution chosen as the one closest to the \textit{n}-th solution. The reachability plot of Extended Data Fig. \ref{fig:figS2} shows the distances between the $n+1$-st and \textit{n}th solutions. Valleys in the reachability plot clearly identify blocks of closely spaced solutions, constituting our solution clans. Clans (valleys) were selected by setting a minimum clan size to $200$ and an upper threshold $0.05$ for the distance between consecutive solutions in a clan. These clans are the clans marked in yellow and blue in Fig. \ref{fig:fig3}b and Extended Data Fig. \ref{fig:figS1}. As the initial solution is chosen at random, this gives \textit{N} distinct ways of sorting the solutions. We found that the large clans were essentially independent of the choice of the initial solution used for the sorting.

\subsection*{Physical interpretation of solution strategies}
In the tunneling clan, marked yellow in Figs. \ref{fig:fig2} and \ref{fig:fig3} and Extended Data Fig \ref{fig:fig3}, the atom is collected by \textit{tunneling} it into the tweezer potential, which is placed on the left hand side of the static potential. As illustrated in Fig. ~\ref{fig:fig3}a, around the time $t=0.06$, all of the $\sim 500$ yellow solutions move to a very particular location in space, which we interpret as the position maximizing the tunneling rate. Instead, the depth of the tweezer potential at this position is only weakly correlated with the final fidelity of the solution (Fig~\ref{fig:fig3}a), permitting large variations in its values. This strategy fails for short durations ($T\leq0.22$) as there is no longer time for the atom to tunnel completely between the potential wells. Although seemingly simple, this strategy is at variance with the intuition obtained from similar tunneling-based problems~\cite{Calarco2004,Anderlini2007,Joergensen2014}, which require careful resonant matching of initial and final state energy levels. In contrast, the CHOP solution utilizes a deep potential for the transport tweezer. Combined with precise motional control of the tweezer, this loads the atom directly into a state with high energy spread (see Fig. ~\ref{fig:fig3}c).

In the shovelling clan, marked blue in Figs. \ref{fig:fig2} and \ref{fig:fig3}, the tweezer is moved past the position of the atom. Therefore the overlap of the tweezer and the static potential forms a strong potential gradient that accelerates the atom towards the target (see Fig. \ref{fig:fig3}b). This strategy transfers the atom into a superposition of many different instantaneous energy levels with a large energy spread $\Delta E$ (see Fig. \ref{fig:fig3}c) allowing fast motion in Hilbert space towards the target state. As illustrated in the inset of Extended Data Fig. \ref{fig:figS3}a, this shovelling strategy stays optimal for shorter durations than the tunneling strategy. Below $T=0.20$, however, it cannot be executed due to a physical bound imposed on the speed of the tweezer, and the yellow tunneling clan becomes superior.

\subsection*{Mapping high-dimensional landscape to 2D \label{sec:sysproj}}
The visualization of the control landscape maps all solutions to a 2D surface. The mapping aims to represent distances $D_{j,l}$ between any two solutions. This distance is a Manhattan type distance,
\[
D_{j,l}=\sum_i \int _0^T |u^j_i - u^l_i| \mathrm{d} t,
\]
where $u^{j(l)}_i$ are rescaled tweezer position and amplitude to intervals of unit length. The solutions are mapped to the 2D-landscape by the construction rules:
\begin{enumerate}
\item A randomly chosen solution is placed at the origin $(x_1,y_1)=(0,0)$
\item The two solutions closest to the initial solution are given Euclidean coordinates $(x_i,y_i)$ such that Euclidean distances between them match the Manhattan type distances $D_{1,2}$, $D_{1,3}$, and $D_{2,3}$. The points define a triangle with an arbitrary orientation.
\item The solution with the smallest Manhattan distance to the previous ones is given coordinates $(x_k,y_k)$ such that the 2D distances are as close as possible to the values  $D_{i,k}$ $i=1,2,..,k-1$. This is done with the Nelder-Mead algorithm using the cost function,
\begin{equation}
S_k=\sum_{j=1}^{k-1}\Bigl|D^E_{j,k}-D_{j,k}\Bigr|,
\end{equation}
where $D^E_{j,k}$ denotes the Euclidean distance between the coordinates of two solutions in 2D landscape.
\end{enumerate}
A 3D figure is obtained by plotting $(x_i,y_i,F_i)$ where $F_i$ is the fidelity of the \textit{i}-th solution. This procedure gives a qualitative impression of the optimization landscape. Landscapes are constructed from player seeds, random computer seeds and optimized solutions at the durations $T\simeq 0.17$ and $T=0.40$ are shown in Fig. \ref{fig:fig4}. The multitude of spikes appears because gradient based optimization leads to nearly vertical lines in the landscape. This highlights the failure of local optimization algorithms to explore extended parts of the global landscape. A similar landscape for HILO is also shown in Fig. \ref{fig:fig4}.

\subsection*{HILO \label{sec:param}}
HILO, or Heuristically Initialized Local Optimizer, applies the Krotov algorithm to a seed found by using the same heuristic as the best CHOP solutions. The motion of the tweezer in a tunneling (yellow) and shoveling (blue) CHOP solution can be reasonably approximated by three movements at constant horizontal speed and rate of change of the potential depth (see Extended Data Fig. \ref{fig:figS3}a), which connects three points $P_i=(x_i,A_i,t_i)$, $i=1,2,3$. These three lines correspond to the first quick movement to a point $P_1$, from $P_1$ a slow backwards motion to $P_2$ and a final move to $P_3$ inside the target area. This defines an optimization problem with dimension $D=9$. $P_3$ was set equal to the last point in $w_i$. We noticed that the duration of the quick movement ($t_1$) should be as small as possible, reducing the seed space dimension to $D=5$. This dimension could be further reduced to $D=3$ by only changing the position of the tweezer using straight lines and letting the amplitude decay exponentially to $\mathcal{A}=-150$ after $t_1$. Optimizations showed that $\mathcal{A}=-150$ was the optimal value.
To reach the desired state with high fidelity snaking the tweezer is crucial. The KA is excellent at introducing the snaking to stabilize the wave function so this parametrization deliberately does not include these. We generated seeds for KA from this space using a direct search. 

\subsection*{The Direct Hilbert Velocity}
In numerous studies the fidelity below the QSL has been found to follow a universal $\sin^2\bigl(\frac{\pi}{2} T/T_{\mathrm{QSL}}\bigr)$ behavior~\cite{Caneva2011}.
This section briefly summarizes how an expression for the fidelity can be obtained using the concept of direct Hilbert velocity \textit{Q} \cite{Gajdacz2014}. State $|\psi(0)\rangle$ is evolved in time under a controlled Hamiltonian $\hat{\mathcal{H}(t)}$ towards a target state $|\chi(T)\rangle$ at time \textit{T}. The Krotov algorithm exploits the fact that the target state can be propagated backwards in time using the (adjoint) time evolution operator. Using $|\chi(t)\rangle$, the process fidelity can be evaluated at any instant of time,
\[
F= |\langle \chi(T) | \psi(T) \rangle|^2 = |\langle \chi(t) | \psi(t) \rangle|^2.
\]
It is useful to define the component $|\xi \rangle$ of the backward evolved target state $|\chi(t)\rangle$, which is orthogonal to the instantaneous state $|\psi(t)\rangle$,
\[
|\xi \rangle = \frac{|\chi\rangle \langle \chi | - F}{\sqrt{F(1-F)}}|\psi \rangle.
\]
The norm of this component reduces at a rate
\[
Q(t)=\mathrm{Re}\Bigl\langle \xi(t)|\dot{\psi}(t)\Bigr\rangle=\mathrm{Im}\Bigl\langle \xi(t)\Bigl|\hat{\mathcal{H}}(t)\Bigr|\psi(t)\Bigr\rangle,
\]
which is denoted the direct Hilbert velocity~\cite{Gajdacz2014}. As described in ~\cite{Gajdacz2014}, the trade-off between process duration and achievable fidelity obeys the relation
\begin{equation}
\frac{\mathrm{d} F}{\mathrm{d} T}=\frac{2 \sqrt{F(1-F)}}{T} \int_0^T Q(t) \mathrm{d}t = 2\sqrt{F(1-F)} \langle Q \rangle_T, \label{QT}
\end{equation}
where $\langle Q \rangle_T$ is the time average of \textit{Q} over the process duration. Equation (\ref{QT}) is always true for a uniform extension of time and it is valid for any extension of time for an optimal process \cite{Gajdacz2014}. Integrating equation (\ref{QT}) for an optimal clan with vanishing initial fidelity leads to
\begin{equation}
F=\sin^2\biggl(\int_0^T \langle Q \rangle_{T'} \mathrm{d} T'\biggr). \label{Qtofidelity}
\end{equation}
Implying that whenever $\langle Q \rangle_T$ is independent of time, $T_{\mathrm{QSL}}=\pi/2\langle Q \rangle_T$ and for
shorter durations, $F(T) = \sin^2\bigl(\frac{\pi}{2} T/T_{\mathrm{QSL}}\bigr)$. However, \textit{ a priori} there is no reason that the average direct Hilbert speed should be constant over an extended range of process durations. $\langle Q \rangle_T$ has been calculated for BHW for the optimal CHOP paths using the method described above, shown in Extended Data Fig. \ref{fig:figS3}b.
Note that $|\xi\rangle$ and hence $\langle Q \rangle_T$ is only defined for durations where $F<1$. Extended Data Figure \ref{fig:figS3}b shows that $\langle Q \rangle_T$ is not constant in time, which explains the deviations from $\sin^2\bigl(\frac{\pi}{2} T/T_{\mathrm{QSL}}\bigr)$ behaviour in Extended Data Fig. \ref{fig:figS3}a. How the fidelity drops when approaching the QSL for the different strategies can be seen in the insert in Extended Data Fig. \ref{fig:figS3}a.

\bibliographystyle{naturemag}
\renewcommand\refname{Additional References and Notes}

\section*{Extended Data Figures}
\setcounter{figure}{0}
\renewcommand{\figurename}{\textbf{Extended Data figure}}

\begin{figure}[H]
\centering
\includegraphics[width=0.6\textwidth]{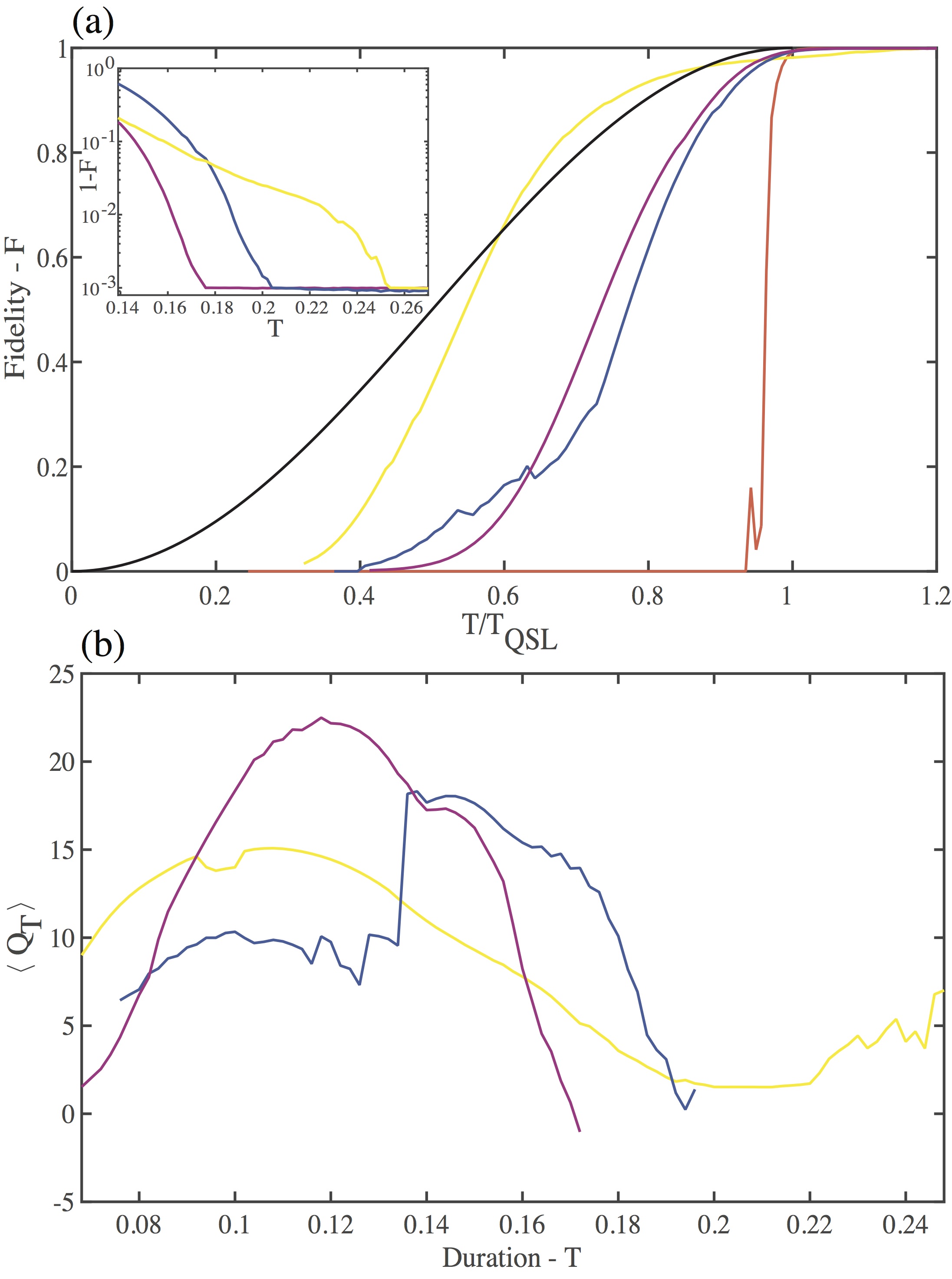}
\caption{\textbf{Deviations from }$\mathbf{sin^2}$\textbf{ - behaviour.} Colouring as in Fig. 1. Panel \textbf{a} shows the fidelity as a function of duration for the optimal families rescaled by the apparent QSL. The $T_{\mathrm{QSL}}$ is found for each solution by fitting $\sin^2(aT+b)$. The black line shows $\sin^2(\pi/2 T)$ for reference. Panel \textbf{b} shows the direct Hilbert speed $\langle Q \rangle_T$ for the best solutions found by HILO, \textit{tunneling} and \textit{shoveling} clans respectively. Note that $\langle Q \rangle_T$ is only defined for durations with $F<1$ so the curves end at different durations. The varying direct Hilbert speed explains the deviation from $\sin^2\bigl(\frac{\pi}{2} T/T_{\mathrm{QSL}}\bigr)$ on the \textbf{a}. }
\label{fig:figS1}
\end{figure}

\begin{figure}[H]
\centering
\includegraphics[width=0.6\textwidth]{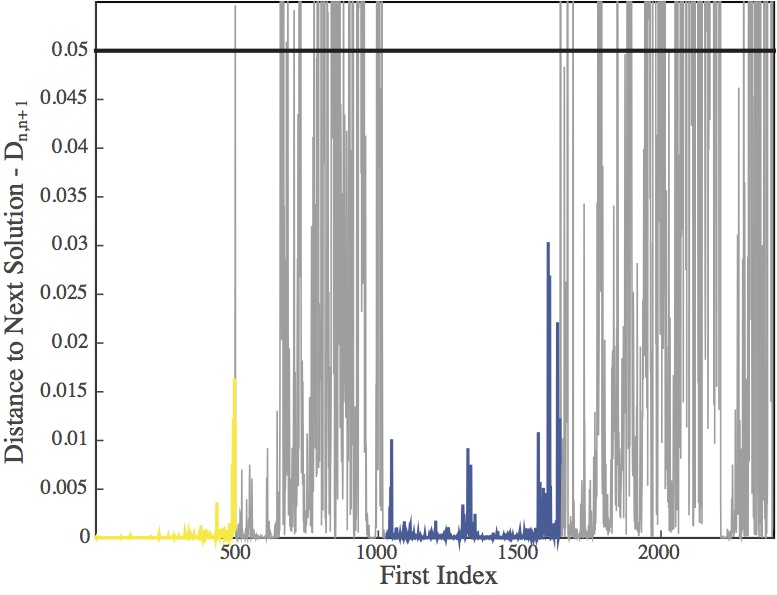}
\caption{\textbf{The reachability plot for Fig. \ref{fig:fig3}b.} The distance between subsequent solutions as calculated by Eq. \ref{distance}. Valleys identify closely spaced solutions, denoted clans. The valleys corresponding to tunneling and shoveling clans are marked with yellow and blue respectively. The black line marks the threshold for the distance between consecutive solutions in a clan at $0.05$.}
\label{fig:figS2}
\end{figure}

\begin{figure}[H]
\centering
\includegraphics[width=0.6\textwidth]{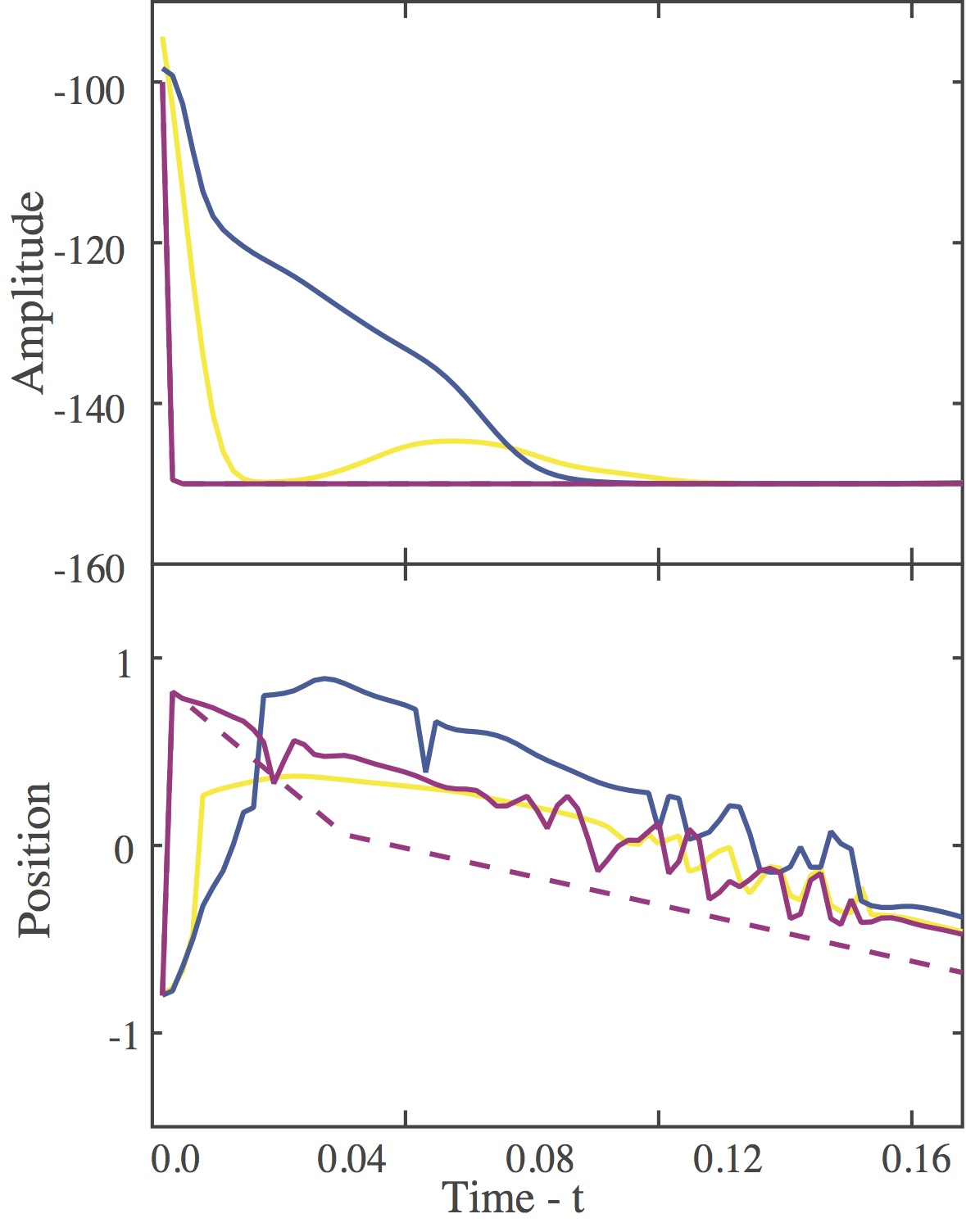}
\caption{\textbf{Solutions from CHOP and HILO.} The amplitude and the position of the tweezer as a function of time for the best player solutions in the tunneling (yellow) and shoveling (blue) clans and the best HILO (purple). The dashed purple line shows the initial seed used by the best HILO solution (note that the dashed and solid purple lines in the top panel overlap). The total duration is $T = 0.15.$}
\label{fig:figS3}
\end{figure}

\end{document}